%% file: doc.tex
\documentclass[11pt]{article}

\usepackage[a4paper, top=50pt, bottom=50pt]{geometry}

\usepackage[utf8]{inputenc}
\usepackage{hyperref, amssymb, amsmath, xcolor, mathtools, makecell, tikz,
            multicol, microtype, xspace, mathabx, shuffle, float}
\usepackage{enumitem}
\usepackage[lofdepth,lotdepth]{subfig}
\usetikzlibrary{matrix, calc, arrows, petri}
\usepackage{racc}
\usepackage{caption} 
\usepackage[thmmarks,amsmath]{ntheorem}
\theoremstyle{plain}
\theoremseparator{:}
\newtheorem{defi}{Definition}[section]
\newtheorem{thm}[defi]{Theorem}
\newtheorem{lemme}[defi]{Lemma}
\newtheorem{coro}[defi]{Corollary}
\newtheorem{prop}[defi]{Proposition}

\theorembodyfont{\upshape}
\newtheorem{ex}{Example}
\newtheorem{remarq}[defi]{Remark}

\newcommand{\longversion}[1]{}

\theoremheaderfont{\sc}
\theoremstyle{nonumberplain} 
\theoremsymbol{\rule{1ex}{1ex}} 
\newtheorem{preuve}{Proof}
\newcommand{\email}[1]{\url{#1}}

\author{Jean-Baptiste Priez\thanks{\email{jean-baptiste.priez@lri.fr}}
\and Aladin Virmaux\thanks{\email{aladin.virmaux@lri.fr}}}

\title{Non-commutative Frobenius characteristic of generalized
parking functions\\Application to enumeration}

\begin{document}
\maketitle
\begin{abstract}\noindent
  \begin{description}
    \item[English] We give a recursive definition of generalized parking
      functions that allows them to be viewed as a species. From there we
      compute a non-commutative characteristic of the generalized parking
      function module and deduce some enumeration formulas of structures and
      isomorphism types. We give as well an interpretation in several bases of
      non commutative symmetric functions. Finally, we investigate an
      inclusion-exclusion formula given by Kung and Yan.
    \item[French] Nous donnons une définition récursive des fonctions de
      parking généralisées nous permettant de munir ces dernières d'une
      structure d'espèce. Nous utilisons ce point de vu pour donner une
      caractéristique de Frobenius non-commutative du module des fonctions de
      parking généralisées que nous appliquons afin de donner de nombreuses
      formules d'énumération de structures et de type d'isomorphismes, ainsi
      qu'une interprétation dans plusieurs bases des fonctions symétriques non
      commutatives. Enfin, nous étudions une formule d'inclusion-exclusion
      provenant de Kung et Yan.
  \end{description}
\end{abstract}
\tableofcontents
\section*{Introduction} 
\label{sec:intro}
\input{intro}
\section{Species} 
\label{sec:background}
\input{background}
\section{Parking functions}
\label{sec:pf}
\input{pf}
\section{Non-commutative Frobenius characteristic of $\PF(\seqPF)$}
\label{sec:ncch}
\input{ncsf}
\section{An inclusion-exclusion formula}
\label{sec:incl-excl}
\input{incl-excl}
\section*{Acknowledgments}
\label{sec:acknowledgement}
\input{acknowledgement}
%
\setlength{\parindent}{-15pt}
\bibliographystyle{alpha}
\bibliography{biblio}
\label{sec:biblio}

\end{document}

%% file: intro.tex
Parking functions were introduced in \cite{konheim1966occupancy} to
model hashing problems in computer science and appear in many
different contexts in combinatorics.
Generalized parking functions were introduced by Stanley and
Pitman~\cite{stanley2002polytope}. Shortly after, Kung
and Yan showed that the Gon$\Breve{\text{c}}$arov polynomials form a
natural basis to manipulate generalized parking functions and stated
numerous enumeration formulas~\cite{kung2003goncarov}.

The (linear span of the) set $\PF_n$ of parking functions of length
$n$ is naturally a module over the symmetric group $S_n$ acting on
positions. In~\cite{novelli2008noncommutative} Novelli and Thibon
observed that $\PF_n$ is also a module over the $0$-Hecke algebra
$H_n(0)$. This allows us to lift canonically the Frobenius characteristic
of $\PF_n$, which lives in the algebra of symmetric functions, as a
non commutative Frobenius characteristic in the algebra
of non commutative symmetric functions. They then apply this
technology to the non commutative Lagrange inversion. Implicit in the
calculations of~\cite{novelli2008noncommutative} is the use of a
recursive definition of parking functions to derive functional
equations on the Frobenius characteristic.

In this paper we apply the same representation-theoretic approach to
generalized parking functions and derive enumeration formulas such as
those of Kung and Yan. A key observation is that generalized parking
functions naturally form a species and that this species can be
defined recursively.

This paper is organized as follows, beginning with background on species in
Section~\ref{sec:background}, we the define in Section~\ref{sec:pf} the species
$\PF(\seqPF)$ of generalized parking functions. Then we will show that
generalized parking functions can be defined recursively, which naturally
translates into a functional equation on $\PF(\seqPF)$
(Theorem~\ref{eq:pf_recursive_definition}). We derive a closed-form expression
for this species by expressing it in terms of the species $\speciesSet$ of
sets.


In Section~\ref{sec:ncch} we apply the previous results to the
computation of the non-commutative Frobenius characteristic of the
module of parking functions, expressed in the complete basis of
non-commutative symmetric functions
(Theorem~\ref{thm:PF_characteristic_expand}). From there, we derive
new enumeration formulas, then 
express $\ncch(\PF(\seqPF))$ in the ribbon and $\Lambda$ basis; the
latter admits a nice combinatorial interpretation
(Proposition~\ref{prop:lambda-combi}).

Finally in Section~\ref{sec:incl-excl} we state and prove an inclusion
exclusion formula on the faces of a polytope for generalized parking functions,
giving a combinatorial interpretation of ~\cite[Theorem 4.2]{kung2003goncarov}.

%% file: background.tex
In this paper we use the theory of species~\cite{bergeron1998combinatorial} to
encode the notion of labeled and unlabeled parking functions simultaneously. We
recall some definitions and some classical operations for species.
\paragraph{Definition:}
A \emph{species} $\speciesA$ is an endofunctor of the category of sets
with bijections into itself.\\\noindent
In other terms, $\speciesA$ is a rule which produces 
\begin{multicols}{2}
\begin{itemize}
  \item a finite set $\speciesA[U]$, for any finite set $U$,
  \item a function $\speciesA[\sigma]$, for any bijection $\sigma : U \to V$.
\end{itemize}
\end{multicols}
This function $\speciesA[\sigma]$ satisfies the functoriality conditions:
\begin{align*}
\speciesA[\tau \circ \sigma] = \speciesA[\tau] \circ \speciesA[\sigma]\,; &&
\speciesA[Id_U] = Id_{\speciesA[U]}\,,
\end{align*}
for any bijections $\sigma : U \to V$ and $\tau : V \to W$,
and with the identity map $Id_U : U \to U$.

Elements of $\speciesA[U]$ are called the \emph{$\speciesA$-structures on $U$}
and functions $\speciesA[\sigma]$ are called the \emph{transports} of
$\speciesA$-structures along $\sigma$.
Two structures $f \in \speciesA[U]$ and $g \in \speciesA[V]$ have the same
\emph{isomorphism type} if there is a bijection $\sigma : U \to V$ such that
$\speciesA[\sigma](f) = g$.
We denote $\speciesA_n$ the species $\speciesA$ \emph{restricted to
sets of cardinality $n$}.
\paragraph{Characteristic species} 
The species $\speciesOne$ is $\speciesOne[U] = \{\emptyset\}$ if $U = \emptyset$
and $\emptyset$ in otherwise.
\paragraph{Species of sets}
The species of sets $\speciesSet$ is defined by $\speciesSet[U] = \{U\}$, for
any finite set $U$ (endowed with the trivial action $U \mapsto \sigma(U)$ for
any bijection $\sigma: U \to V$ and any finite set $U$).

Many operations on species allows a direct translation in terms of
generating series: \emph{addition}, \emph{multiplication}, \emph{substitution},
\textit{etc}. These operations constitute combinatorial analogs of the usual
operations on series. Here we will only use \emph{addition} and
\emph{multiplication}. 
In the sequel, let $\speciesA$ and $\speciesB$ be two
species.
\paragraph{Addition}
The sum of $\speciesA$ and $\speciesB$, noted $\speciesA + \speciesB$, is defined
by:
\begin{align*}
    (\speciesA + \speciesB) [U] = \speciesA[U] \sqcup \speciesB[U]\,; &&
    (\speciesA + \speciesB)[\sigma](f) = \begin{dcases*}
        \speciesA[\sigma](f) & if $f \in \speciesA[U]$,\\
        \speciesB[\sigma](f) & if $f \in \speciesB[U]$,
    \end{dcases*}
\end{align*}
for any finite set $U$, any bijection $\sigma : U \to V$ and any $f \in
(\speciesA + \speciesB)[U]$. 
\paragraph{Product}
Similarly the \emph{product of species} of $\speciesA$ and $\speciesB$, noted
$\speciesA\cdot \speciesB$ defines ordered pairs of structures $f = (g, h)$:
\begin{align*}
    (\speciesA \cdot \speciesB) [U] = \sum_{S \sqcup T = U} \speciesA[S] \times
    \speciesB[T]&& (\speciesA \cdot \speciesB)[\sigma](f) =
    (\speciesA[\sigma_g](g), \speciesB[\sigma_h](h))\,,
\end{align*}
with $\sigma_g$ (respectively $\sigma_h$) the restriction of $\sigma$ to the
underlying set of the $\speciesA$-structure $g$ (\resp $\speciesB$-structure
$h$). 
We denote $\speciesA^n$ the product of $\speciesA$ with itself $n$ times:
$\speciesA^n = \speciesA \cdots \speciesA^{n-1}$ (with $\speciesA^1 =
\speciesA$ and $\speciesA^0 = \speciesOne$).

%% file: pf.tex
A \emph{parking function} on finite set $U$ (of cardinality $u$) is a function
$f: U \to \NN_+$ such that $\card{f^{-1}([k])} \seq k$, for any $k \in [u]$. The
generalization introduced in \cite{stanley2002polytope} modifies the condition.
Let $\seqPF : \NN_+ \to \NN$ be a non-decreasing sequence; the function $f
: U \to \NN$ is a \emph{$\seqPF$-parking function} if
$\card{f^{-1}([\seqPFn{k}])} \seq k$, for any $k \in [u]$.
\begin{remarq}
    The usual parking functions are $\seqPF$-parking functions with $\seqPF : i
    \mapsto i$, the identity map.
\end{remarq}
For the following, it will be more convenient to use another equivalent
definition. A (generalized) parking function $f : U \to \NN_+$ may be
described as an ordered sequence of sets $(Q_i)_{i \in \NN_+}$ where $Q_i =
f^{-1}(i)$. From this definition we remark that the 
isomorphism types (unlabeled structures) of generalized parking functions are
(generalized) Dyck paths, namely a staircase walk under a discrete curve
$\seqPF$.
Formally this defines a species:
\begin{defi}
    The species of $\seqPF$-parking functions is: 
    \begin{itemize}
      \item for any finite set $U$, the set of all sequences $(Q_i)_{i \in
      \NN_+}$ of disjoint subsets of $U$ such that
\begin{align}
    \sum_{i = 1}^{\seqPFn{k}} \card{Q_i} \seq k && \text{for any} \quad 1 \ieq
    k \ieq u\,,
    \label{eq:pf_set_condition}
\end{align}
      \item for any bijection $\sigma : U \to V$, the relabeling action is
      $(Q_i)_i \mapsto (\sigma(Q_i))_i$ (for any $\seqPF$-parking function on
      $U$).
    \end{itemize}
\end{defi}
This viewpoint on generalized parking functions in terms of sequences
reveals a recursive definition.
A simple way to put forward the recurrence is to view a generalized parking
function $(Q_i)$ as a decorated path/staircase walk defined by: the first
\emph{tread} (horizontal step) goes to $(\card{Q_1}, 1)$ and is decorated by
$Q_1$. The second tread decorated by $Q_2$ starts at $(\card{Q_1}, 1)$ and goes
to $(\card{(Q_1 \cup Q_2)}, 2)$. The third starts where the second ends, goes
to $(\card{(Q_1 \cup Q_2 \cup Q_3)}, 3)$ and is decorated by $Q_3$,
\textit{etc.} (see Example \ref{ex:pf_representation}). 
\begin{figure}[h]
\begin{ex}
    \label{ex:pf_representation}
    Let $\seqPF$ be the sequence $2,2,3,5,8,8,8,8,8,8,9,9,9,9,9,9,10,
    10, \ldots$ and it is pictured in red. Let $(Q_i) = (\{d,f,m\}, \{g\},
    \{a,h\}, \emptyset, \{c,e,j,k,n,p\},\emptyset, \emptyset, \emptyset,
    \{b,i,l,o,q,r,s,t\}, \emptyset)$ be a
    $\seqPF$-parking function on $\{a, b, c,\cdots, t\}$. The parking function
    is represented by the decorated blue path.
\begin{align*}
    \vcenter{\hbox{\begin{tikzpicture}[scale=.25]
   \def\x{20}
   \def\y{20}
    \foreach \i in {0,...,\x} { \draw [very thin,gray] (\i,0) -- (\i,\y); }
    \draw [very thin,gray] (0,0) -- (\x,0);
    \foreach \i in {2,4,...,\y} {
        \draw [very thin,gray] (0,\i) -- (\x,\i);
    }
    \node[red] at (3, 14) {$\seqPF$};
    \draw[thick, red] (0,0) -- (0, 4) -- (2,4) -- (2,6) -- (3,6) -- (3,10) 
                            -- (4, 10) -- (4, 16) -- (10, 16) -- (10, 18) 
                            -- (16, 18) -- (16, 20) -- (20, 20);
    \draw[thick, blue!60] (0, 2) -- (3,2) -- (3,4) -- (4,4) -- (4,6) -- (6, 6)
    -- (6, 10) -- (12, 10) -- (12,17.9) -- (20, 17.9) -- (20,20); 
    \node at (1.7,1) {$\scriptstyle \{d,f,m\}$}; 
    \node at (3.7,3) {$\scriptstyle \{g\}$};
    \node at (5.2,5) {$\scriptstyle \{a,h\}$};
    \node at (6.7,7) {$\scriptstyle \emptyset$};
    \node at (9.2,9) {$\scriptstyle \{c,e,j,k,n,p\}$};
    \foreach \i in {11,13,15}
        \node at (12.7,\i) {$\scriptstyle \emptyset$};
    \node at (16,17) {$\scriptstyle \{b,i,l,o,q,r,s,t\}$};
    \node at (20.5,19) {$\scriptstyle \emptyset$};
    \end{tikzpicture}}}
\end{align*}
\end{ex}
\end{figure}

From this graphic representation it is easy to notice that a
$\seqPF$-parking function $(Q_i)$ on $U$ is either a sequence of empty sets if
$U = \emptyset$, or a sequence of sets $(Q_i)_{i \in [\seqPFn{1}]}$ of union
$S\subseteq U$, concatenated with another generalized parking function on $U
\backslash S$ defined from another non-decreasing function $\seqPFb_s$ (See
Example \ref{ex:split_pf}). This new map is defined by a \emph{shift} of
$\seqPF$ characterized by the cardinality $s$ of $S$:
\begin{align*}
    \label{eq:seqPFb_en_fonction_de_seqPF_n}
    \seqPFbn[s]{m} = \seqPFn{s+m} - \seqPFn{1}\,.
\end{align*}
In the following, we will implicitly denote by $\seqPFb_n$ the shift of
$\seqPF$ by $n$.
\begin{figure}[h]
\begin{ex}
    \label{ex:split_pf}
    Let $(Q_i)$ be the $\seqPF$-parking function defined in Example
    \ref{ex:pf_representation}; it is the concatenation of
    the sequence of the two first sets $\boldsymbol{(}\{d,f,m\},
    \{g\}\boldsymbol{)}$ and the $\seqPFb$-parking function $\boldsymbol{(}\{a,h\}, \emptyset,
    \{c,e,j,k,n,p\},\emptyset, \emptyset, \emptyset, \{b,i,l,o,q,r,s,t\},
    \emptyset\boldsymbol{)}$ on $\{a,b,c, \cdots, t\} \backslash \{d,f,g,m\}$
    with $\seqPFb$ the sequence $3,6,6,6,6,6,6,7,7,7,7,7,7, 8, 8, \ldots$
    In other terms we have $\seqPFbn{m} = \seqPFn{m + 4} - \seqPFn{1}$.
    \begin{align*}
    \vcenter{\hbox{\scalebox{.8}{\begin{tikzpicture}[scale=.3]
   \def\x{20}
   \def\y{20}
    \foreach \i in {0,...,\x} { \draw [very thin,gray] (\i,0) -- (\i,\y); }
    \draw [very thin,gray] (0,0) -- (\x,0);
    \foreach \i in {2,4,...,\y} {
        \draw [very thin,gray] (0,\i) -- (\x,\i);
    }
    \node[red] at (3, 14) {$\seqPF$};
    \draw[thick, red] (0,0) -- (0, 4) -- (2,4) -- (2,6) -- (3,6) -- (3,10) 
                            -- (4, 10) -- (4, 16) -- (10, 16) -- (10, 18) 
                            -- (16, 18) -- (16, 20) -- (20, 20);
    \draw[thick, blue!60] (0, 2) -- (3,2) -- (3,4) -- (4,4) -- (4,6) -- (6, 6)
    -- (6, 10) -- (12, 10) -- (12,17.9) -- (20, 17.9) -- (20,20); 
    \node at (1.7,1) {$\scriptstyle \{d,\ f,\ m\}$}; 
    \node at (3.7,3) {$\scriptstyle \{g\}$};
    \node at (5.2,5) {$\scriptstyle \{a,\ h\}$};
    \node at (6.7,7) {$\scriptstyle \emptyset$};
    \node at (9.2,9) {$\scriptstyle \{c,\ e,\ j,\ k,\ n,\ p\}$};
    \foreach \i in {11,13,15}
        \node at (12.7,\i) {$\scriptstyle \emptyset$};
    \node at (16,17) {$\scriptstyle \{b,\ i,\ l,\ o,\ q,\ r,\ s,\ t\}$};
    \node at (20.5,19) {$\scriptstyle \emptyset$};
    \end{tikzpicture}}}} = 
    \boldsymbol{(}\{d,f,m\}, \{g\}\boldsymbol{)} \cdot
    \vcenter{\hbox{\scalebox{.8}{\begin{tikzpicture}[scale=.3]
   \def\x{20}
   \def\y{20}
    \foreach \i in {4,...,\x} { \draw [very thin,gray] (\i,4) -- (\i,\y); }
    \foreach \i in {4,6,...,\y} {
        \draw [very thin,gray] (4,\i) -- (\x,\i);
    }
    \node[red] at (3, 14) {$\seqPFb$};
    \draw[thick, red] (4, 4) -- (4, 16) -- (10, 16) -- (10, 18) 
                            -- (16, 18) -- (16, 20) -- (20, 20);
    \draw[thick, blue!60] (4,6) -- (6, 6)
    -- (6, 10) -- (12, 10) -- (12,17.9) -- (20, 17.9) -- (20,20); 
    \node at (5.2,5) {$\scriptstyle \{a,\ h\}$};
    \node at (6.7,7) {$\scriptstyle \emptyset$};
    \node at (9.2,9) {$\scriptstyle \{c,\ e,\ j,\ k,\ n,\ p\}$};
    \foreach \i in {11,13,15}
        \node at (12.7,\i) {$\scriptstyle \emptyset$};
    \node at (16,17) {$\scriptstyle \{b,\ i,\ l,\ o,\ q,\ r,\ s,\ t\}$};
    \node at (20.5,19) {$\scriptstyle \emptyset$};
    \end{tikzpicture}}}}
\end{align*}
\end{ex}
\end{figure}
\subsection{Recursive definition}
The recursive splitting described earlier involves a natural constructive
definition of the \emph{$\seqPF$-parking functions} species in terms of species
operations based on the species of sets $\speciesSet$.
\begin{thm}
    \label{thm:PF_recursive_construction}
    The species of $\seqPF$-parking functions is isomorphic to the
    species $\PF(\seqPF)$ recursively defined as
    \begin{align}
    \PF(\seqPF) &= (\speciesSet^{\seqPFn{1}})_0 + 
                   \sum_{n \seq 1} (\speciesSet^{\seqPFn{1}})_n \cdot
                   \PF(\seqPFb_n).
          \label{eq:pf_recursive_definition}
    \end{align}
\end{thm}
\begin{preuve}
    A generalized parking function $(Q_i)$ is an infinite sequence of disjoint
    subsets of a finite set $U$ and with $Q_i = \emptyset$ for $i >
    \seqPFn{u}$, so $(Q_i)$ may be seen as a sequence of length $\seqPFn{u+1}$.
    The relabeling action is trivially the same.
    By induction, any $\PF(\seqPF)$-structures satisfies the
    generalized parking functions condition (\ref{eq:pf_set_condition}).
    Finally, any $\seqPF$-parking functions on $U$ can be divided
    into factors $D_p = (Q_i)_{i}$ with $\seqPF(p-1) < i \ieq \seqPF(p)$ for
    any $p \in [u]$. By induction again, each factor $D_p$ corresponds to a
    structure on the left term of the sum \ref{eq:pf_recursive_definition},
    ($(\speciesSet^{\seqPFbn[{\alpha_p}]{1}})_{\beta_p}$) with
    $\alpha_p = \beta_{p-1}$.
\end{preuve}
Thanks to species theory \cite{bergeron1998combinatorial}, this constructive
definition of the $\seqPF$-parking functions gives automatically a cycle index
series: the series of (commutative) Frobenius characteristic of the
natural symmetric group action on $\PF_n(\seqPF)$.
Furthermore, the terminal elements of our grammar $(\speciesSet^k)_n$ are
well-known to be characterized by the \emph{permutational $0$-Hecke modules}.
\subsection{Closed-form equivalent definition}
In order to give a closed-form expression of the species $\PF(\seqPF)$ one
needs to understand the map $\seqPFb_n$ according to $\seqPF$. Tracking the
recursion of the definition is only about following the different choices of $n$
in
(\ref{eq:pf_recursive_definition}).
This recursive definition is mainly defined by the map $\seqPFb_n$ of equation
(\ref{eq:seqPFb_en_fonction_de_seqPF_n}). When expanding $\PF(\seqPF)$, we
remark that the map evolves as follows:
\begin{align*}
    \seqPFb_n : m \quad\mapsto \quad&\seqPF(m+n) - \seqPF(1)\\
    \seqPFb_{n_1}': m \quad\mapsto \quad& \seqPF(m+n+n_1) - \seqPF(1+ n_1)\\
    \seqPFb_{n_2}'': m \quad\mapsto \quad& \seqPFn{m+ n + n_1 + n_2} - \seqPFn{1
    + n_1 + n_2}\\
    \vdots
\end{align*}
On the other hand each $\seqPFb_n$ is always called with $m=1$ by definition and
each $n$ is always a positive integer. So the sequence $(n, n_1, n_2, \cdots)$
is finite sequence of positive integer: a composition $\pi$.
\begin{defi}
  Let $\pi = \pi_1 \cdots \pi_k$ be a composition of $n$. The map $\Upsilon$
  is defined by:
\begin{align}
    \Upsilon(\seqPF; \pi, i) = \begin{dcases*}
        \seqPFn{1} & if $i = 1$,\\
        \seqPFn{1 + \pi(i-1)} - \seqPFn{1 + \pi(i-2)} & otherwise. 
    \end{dcases*}
    \label{eq:upsilon_def}
\end{align}
with $\pi(i) = \pi_1 + \cdots + \pi_i$ the partial sum of the first $i$ parts of
$\pi$.
\end{defi}
By expanding the recurrence (\ref{eq:pf_recursive_definition}) of
theorem~\ref{thm:PF_recursive_construction} we have:
\begin{prop}
    \label{prop:pf_set_def}
    \begin{align}
    \PF(\seqPF) &= \speciesOne + \sum_{n \seq 1} \PF_n(\seqPF) &
\text{with}\qquad
    \PF_n(\seqPF) = \sum_{\pi \models n} \prod\limits_{i=1}^{\ell(\pi)}
                    \left( \speciesSet^{\Upsilon(\seqPF; \pi, i)}
                    \right)_{\pi_i}\,.
    \label{eq:cor_pf_set_def}
    \end{align}
\end{prop}
\begin{preuve}
    The exact formula obtained is
    \begin{align*}
    \PF(\seqPF) = (\speciesSet^{\seqPFn{1}})_0 + 
    \sum_{n \seq 1}\left(\sum_{\pi \models n} \prod_{i=1}^{\ell(\pi)} \left(
        \speciesSet^{\Upsilon(\seqPF; \pi, i)} \right)_{\pi_i} \right) \cdot 
                    \left(\speciesSet^{\seqPFn{n+1} - \seqPFn{n}}\right)_0.
    \end{align*} 
    It can be simplified by turning all pending empty sets at the end into the
    species $\speciesOne$ (remember that $\speciesOne$ is a neutral element for
    product of species).
\end{preuve}

%% file: ncsf.tex
In species theory, there are many combinatorial operations on structures, which
are translated on operations in the cycles index series. In the case of
generalized parking functions, the functional equation/grammar
(\ref{eq:pf_recursive_definition}) is terminating on (and only on) finite
sequences of sets. Those structures are well-known to have a more expressive
non-commutative Frobenius characteristic.

We recall those characteristic in the first subsection, we then give the
non-commutative characteristic of $\PF(\seqPF)$ in bases: $(\S^\pi)$ the
completes, $(\R_\pi)$ the ribbons Schur and $(\L^\pi)$ the elementaries of the
non-commutative symmetric functions. (Refer to
\cite{gelfand1994noncommutative} for an overview on non-commutative symmetric
functions.)
\subsection{Species of sequence of $k$-sets}
\label{ssec:species_of_sequences_ksets}
In this subsection we focus on the species $\left(\speciesSet^{\Upsilon(\seqPF;
\pi, i)}\right)_{\pi_i} \simeq \left(\speciesSet^k\right)_n$. 
In \cite{krob1997noncommutative}, the authors lift the right
action of $\mathfrak{S}_n$ on $[k]^n$ by considering the natural right action of
$H_n(0)$ on $\CC[k]^n$. In the same way, we consider here the
natural action of $H_n(0)$ on the linearized species $\CC\speciesSet^k[n]$.
Using species theory notations, we translate some classical results
appearing in~\cite{krob1997noncommutative} (and~\cite{novelli2008noncommutative}). 

Let $Q = (Q_i)_{i \in [k]}$ be a structure in $\speciesSet^k[n]$ that
is a sequence of $k$ disjoint subsets which covers the finite set $[n]$; more
generally we could replace $[n]$ by any finite set $U$ endowed with a
fixed total order so that the elementary transpositions are well
defined.
The Hecke algebra $H_n(0)$ acts on $\CC \speciesSet^k[n]$ on the left
by permuting the elements.  By abuse of notations we note $Q^{-1}(i)
= k$ if $i \in Q_k$. For $q=0$, the action of $T_i$ is defined by:
\begin{align}
  \label{eqn:action-0-hecke}
  Q \cdot T_i =
  \left\{
    \begin{array}{c l}
    \sigma_i(Q)  & \text{if}\ Q^{-1}(i) < Q^{-1}(i+1)\\
    0 & \text{if}\ Q^{-1}(i) = Q^{-1}(i+1)\\
    -Q & \text{otherwise,}
  \end{array}
  \right.
\end{align}
where $\sigma_i$ is the corresponding elementary
transposition ($\sigma_i$ is defined as the bijection with $i \mapsto i+1$,
$i+1 \mapsto i$ and stays fixed otherwise ); for example
$\boldsymbol{(}13|\cdot|2\boldsymbol{)} \cdot T_1 =
\boldsymbol{(}23|\cdot|1\boldsymbol{)}$.
%

The action of any element $T_i$ of $H_n(0)$ on $Q$ is either $0$ or a
rearrangement of $Q$.
The orbits (isomorphism types) are indexed by decompositions $d =
(d_1,\ldots,d_k)$ of $n$ in $k$ parts with $d_i = \card{Q_i}$, that is
generalized compositions including null parts.  The rearrangements of
$Q$ form a basis of an $H_n(0)$-projective module $M$ whose
non-commutative characteristic is $\ncch(M) = S^{\pi}$, where $\pi$ is
the underlying composition of $d$ obtained by stripping away null parts.

The non-commutative characteristic of $\CC\speciesSet^k[n]$ is therefore
\begin{align}
    \ncch(\CC \speciesSet^k[n]) = \sum_{\pi \models n} \mathbf{M}_\pi(k) \S^\pi =
    \sum_{\pi \models n} \binom{k}{\ell(\pi)} \S^\pi = 
    \S_n(k\AA)\,,
    \label{eq:ncch_Q}
\end{align}
where the binomial coefficients account for the number of ways to
insert $k-\ell(\pi)$ empty sets in a sequence of $\ell(\pi)$ non-empty sets.
The non-commutative complete function $\S^\pi$ is used here as a way to encode the
relabeling action of a sequence of $\ell(\pi)$-sets with $\pi_1$ elements
in the first set, $\pi_2$ elements in the second set \textit{etc}. In the enumeration
formula of structures (\ref{eq:enum_struct_E_k}), $\S^\pi(\mathbb{E})$ is specialized
into the multinomial $\binom{n}{\pi_1, \cdots, \pi_j}$ (the reader may consult
the specialization $\S^\pi(\mathbb{E})$ in \cite{hivert20111}).
In terms of Hopf algebras operations $\S_n(k\AA)$ is equivalent to the
Adams operations which iterate $k$ times the coproduct and then the product: $\mu^{k} \circ
\Delta^{k}(\S_n)$ with $\Delta^k = (\Delta \otimes Id^{k-1\otimes}) \circ
\Delta^{k-1}$ and $\mu^k = \mu \circ (Id \otimes \mu^{k-1})$.
\paragraph{$\speciesSet^k$-structures enumeration}
From the characteristic of permutation representations, we recover easily the
enumeration formula of $(\speciesSet^k)_n$-structures (or words on $[k]$ of
length $n$):
\begin{align}
    \structSet_k(n) = \sum_{\substack{\pi \models n\\ \pi = \pi_1 \cdots \pi_j}}
    \binom{k}{j}\binom{n}{\pi_1, \cdots, \pi_j} = k^n\,.
    \label{eq:enum_struct_E_k}
\end{align}
\paragraph{$\speciesSet^k$-isomorphism types enumeration}
Similarly we recover the enumeration formula of
$(\speciesSet^k)_n$-isomorphism types (or non-decreasing words) by specializing
$\S^\pi \mapsto 1$:
\begin{align}
    \tilde{\structSet}_k(n) = \sum_{\pi \models n} \binom{k}{\ell(\pi)} = 
    \binom{n+k-1}{k-1}\,.
    \label{eq:enum_types_E_k}
\end{align}
\subsection{Complete basis formula}
\label{subsec:formulas}
Using~\eqref{eq:ncch_Q} and the recursive definition
\eqref{eq:pf_recursive_definition} we naturally obtain a recursive formula for
the non-commutative Frobenius characteristic series of the $\seqPF$-parking
functions:
\begin{align}
\begin{aligned}
    \ncch(\PF(\seqPF)) &= 1 + \sum_{n \seq 1}
    \ncch(\CC \speciesSet^{\seqPFn{1}}[n]) \ncch(\PF(\seqPFb_n))\nonumber\\
    &= 1 + \sum_{n \seq 1} \S_n(\seqPFn{1}\AA) \ncch(\PF(\seqPFb_n))\,.
\end{aligned}
    \label{eq:characteristic_PF}
\end{align}
By specializing $\S^\pi$ to $\binom{n}{\pi_1, \cdots, \pi_k}$
and to $1$, we obtain (new) formulas to enumerate
$\PF(\seqPF)$-structures and types. Namely from
\eqref{eq:enum_struct_E_k}, we obtain the following recursive
enumeration formula for the number $\structPF(\seqPF; n)$ of
$\PF(\seqPF)$-structures on a set of cardinality $n$:
\begin{align}
    \structPF(\seqPF; n) &
        = \sum_{k = 1}^{n} \binom{n}{k} \seqPFn{1}^k
    \structPF(\seqPFb_k; n-k)
    \label{eq:enum_struct_PF}
\end{align}
with $\structPF(\seqPF; 0) = 1$.  Similarly, we derive from
\eqref{eq:enum_types_E_k} the number of isomorphism types:
\begin{align}
    \typePF(\seqPF; n) &
    = \sum_{k=1}^{n} \binom{n-\seqPF(1)-1}{\seqPF(1)-1} \typePF(\seqPFb_k;
    n-k)\,,
    \label{eq:enum_types_PF}
\end{align}
also with $\typePF(\seqPF; 0) = 1$.
From Proposition \ref{prop:pf_set_def}, we have a non-recursive version of
$\ncch(\PF(\seqPF))$:
\begin{lemme}
    \label{lem:enumeration_formulas}
\begin{align*}
   \ncch(\PF_n(\seqPF)) &= \sum_{
        \pi \models n
   } \ncprod_{i=1}^{\ell(\pi)} \S^{\pi_i}\left(\Upsilon(\seqPF; \pi,i)\AA
   \right)\,.
\end{align*}
\end{lemme}
In \SOld \ref{ssec:species_of_sequences_ksets} we stated that $\S_n(k\AA)$ is
given by the non-commutative Cauchy identity (\ref{eq:ncch_Q}). 
%
This characteristic, expressed as a sum of products of Adams operations
according to $\seqPF$, lifts trivially \cite[Corollary 5.6]{kung2003goncarov}
in non-commutative symmetric functions:
%
\begin{prop}
    \label{prop:lift_Yan_formula_kseq}
    Let $\seqPF, \seqPFb$ be two non-decreasing functions such that $\seqPF(m)
    = \alpha \seqPFb(m)$, for any $m \in \NN_+$.
\begin{align*}
    \ncch(\PF(\seqPF)) &= \ncch(\PF(\seqPFb))(\alpha \AA)\,.   
\end{align*}
\end{prop}
By expanding the
formula of Lemma \ref{lem:enumeration_formulas} we now have a new sum
over compositions, where terms are products of binomials on parts of each
composition (see Table \ref{tab:enum_charact}). 
To get rid of any specialization alphabet, we first need to refine $\Upsilon$
into $\Psi_{\tau}$; namely for $\pi$ a composition of $n$ and $\tau$
a composition of $\ell(\pi)$ we set
\begin{align*}
\Psi_\tau (\seqPF; \pi, i) &= \begin{dcases*}
    \seqPFn{1} & if $i= 1$,\\
    \seqPFn{1 + \pi(\tau(i))} - \seqPFn{1 + \pi(\tau(i-1))} & otherwise\,.
\end{dcases*}
\end{align*}
\begin{remarq}
$\Psi_{(1, \ldots, 1)} = \Upsilon$.
\end{remarq}
The non-commutative characteristic $\ncch(\PF_n(\seqPF))$ of
Lemma~\ref{lem:enumeration_formulas} can now be expanded into the following
theorem:
\begin{thm}
  \label{thm:PF_characteristic_expand}
  The non-commutative characteristic of $\PF_n(\seqPF)$ is given by:
  \begin{align}
    \ncch(\PF_n(\seqPF)) &= \sum_{\pi \models n} \gamma_\pi \S^\pi\,,
  &\text{with}\qquad
    \gamma_\pi = \sum_{\tau \models \ell(\pi)} \prod_{i = 1}^{\ell(\tau)}
    \binom{\Psi_\tau(\seqPF; \pi, i)}{\tau_i}\,.
    \label{eq:gamma_pi}
\end{align} 
\end{thm}
Using again Proposition~\ref{prop:pf_set_def}, from
(\ref{eq:enum_struct_PF}) and (\ref{eq:enum_types_PF}), we get the
following non recursive enumeration formula for $\seqPF$-structures
and isomorphism types:
\begin{align*}
    \structPF(\seqPF; n) = \sum_{
        \substack{\pi \models n\\
        \pi = \pi_1 \cdots \pi_k
   }} \binom{n}{\pi_1, \cdots, \pi_k}\prod_{i=1}^{k}
   \Upsilon(\seqPF; \pi, i)^{\pi_i}\,,&&
   \typePF(\seqPF; n) = \sum_{
        \pi \models n
   } \prod_{i=1}^{\ell(\pi)} \binom{n-\Upsilon(\seqPF; \pi,
   i)-1}{\Upsilon(\seqPF; \pi, i)-1} \,.
\end{align*} 
\begin{ex}
The first values of the non-commutative characteristic of
$\PF(m^2-m+1)$ are given by:
\begin{align*}
    \ncch(\PF_1(m^2-m+1)) &= \S^{1}\,,\qquad
    \ncch(\PF_2(m^2-m+1)) = 2\S^{11} + \S^{2}\\
    \ncch(\PF_3(m^2-m+1)) &= 9\S^{111} + 2\S^{12} + 6\S^{21} + \S^{3}\\
    \ncch(\PF_4(m^2-m+1)) &= 70\S^{1111} + 9\S^{112} + 21\S^{121} + 
                                2\S^{13} + 51\S^{211} + 6\S^{22} + 12\S^{31}
                                + \S^{4}
\end{align*}
\end{ex} 
\begin{table}
\caption{\label{tab:enum_PF_struct}Some enumerations of $\PF(\seqPF)$-structures
for sets of cardinality $n=0$ to $7$.}
\begin{center}
$\begin{array}{c||cccccccc||c}
    \seqPF\backslash n & 0 & 1 & 2 & 3 & 4 & 5 & 6 & 7 & \text{OEIS}\\ \hline
m & 1&  1& 3& 16&  125&  1296& 16807& 262144& \oeis*{A000272}\\
m+1 & 1 & 2 & 8 & 50 & 432& 4802& 65536& 1062882& \oeis*{A089104}\\
2m & 1 & 2 & 12 & 128& 2000& 41472& 1075648& 33554432& 
\oeis*{A097629}\\
%
m^2+m & 1& 2& 20& 512& 25392& 2093472& 260555392& 45819233280& 
\oeis*{A103353}\\
\end{array}$
\end{center}
\end{table}
\begin{table}
\caption{\label{tab:enum_PF_types}Some enumerations of $\PF(\seqPF)$-types for
$n=0$ to $8$.}
\begin{center}
$\begin{array}{c||ccccccccc||c}
    \seqPF\backslash n & 0 & 1 & 2 & 3 & 4 & 5 & 6 & 7 & 8 & \text{OEIS}\\
    \hline 
 m & 1& 1& 2& 5& 14& 42& 132& 429& 1430& \oeis*{A000108}\\
m+2 & 1& 3& 9& 28& 90& 297& 1001& 3432& 11934 &\oeis*{A000245}\\
2m & 1& 2& 7& 30& 143& 728& 3876& 21318& 120175& 
\oeis*{A006013}\\
%
\lceil\frac{m+1}{3}\rceil &   1& 1& 1& 2& 3& 4& 9& 15& 22&
\oeis*{A124753}
\end{array}$
\end{center}
\end{table}
\begin{table}
\caption{\label{tab:enum_charact}The first values of the non-commutative
characteristic with $\seqPFn{m} = a_m + \seqPFn{m-1}$ and $\seqPFn{1} = a_1$ in the complete basis.}
\begin{align*}
    \ncch(\PF_0(\seqPF)) &= 1\,,\qquad
    \ncch(\PF_1(\seqPF)) = a_1\S^1\,,\qquad 
    \ncch(\PF_2(\seqPF)) = a_1\S^2 + \left[ \binom{a_1}{2} + a_1a_2
    \right]\S^{11}\\ 
    \ncch(\PF_3(\seqPF)) &= a_1\S^3 + 
           \left[ \binom{a_1}{2} + a_1(a_2 + a_3)\right] \S^{21} + 
           \left[ \binom{a_1}{2} + a_1a_2\right] \S^{12}\\ &\;\;\;+ 
           \left[ \binom{a_1}{3} + \binom{a_1}{2}(a_2 + a_3) +
           a_1\binom{a_2}{2} + a_1a_2a_3 \right] \S^{111}
\end{align*}
\end{table}
We now investigate how the formula of
Theorem~\ref{thm:PF_characteristic_expand} translates in other natural bases
of non-commutative functions.
%
\subsection{Ribbon Schur basis formula}
%
Recall that the change of basis from the complete basis to the ribbon
Schur functions basis $(\R_\pi)$ is given by $\S^\pi = \sum_{\tau
  \preceq \pi} \R_\tau$, where $\preceq$ denotes the \emph{reverse
  refinement order}.
\begin{ex}
    The compositions $\tau$ of $5$ such that $\tau \preceq 212$ are $212$, $32$,
    $23$ and $5$.
\end{ex}
This change of basis gives the formula:
\begin{align*}
    \ncch(\PF_n(\seqPF)) &= \sum_{\pi \models n} \left(\sum_{\pi \preceq \tau}
        \gamma_\tau
    \right) \R_\pi 
\end{align*}
\begin{ex}
  The first values of the non-commutative characteristic of
  $\PF(2m-1)$ are given by:
    \begin{align*}
\ncch(\PF_1(2m-1)) &= \R_{1}\,,\qquad
\ncch(\PF_2(2m-1)) = 2\R_{1 1} + 3\R_{2}\\
\ncch(\PF_3(2m-1)) &= 5\R_{1 1 1} + 7\R_{1 2} + 9\R_{2 1} + 12\R_{3}\\
\ncch(\PF_4(2m-1)) &= 14\R_{1 1 1 1} + 19\R_{1 1 2} + 23\R_{1 2 1} + 30\R_{1 3} +
28\R_{2 1 1} + 37\R_{2 2} + 43\R_{3 1} + 55\R_{4}
    \end{align*}
The coefficients of $\R_{1^n}$ are (as excepted) Catalan numbers \oeis{A000108},
and the coefficients of $\R_n$ are (less excepted) the number of non-crossing
trees with $n$ nodes \oeis{A001764}.
\end{ex}
\subsection{Lambda basis formula}
Recall that both bases $(\S^\pi)$ and $(\L^\pi)$ are multiplicative
and related by the formula $\S^{\pi} = \sum_{\tau \preceq \pi}
(-1)^{\ell(\tau) - \ell(\pi)} \L^{\tau}$. Furthermore, the change of
base from ribbon $(\R_\pi)$ to lambda $(\L^\pi)$ is given by $\R_\pi =
\sum_{\widebar{\tau} \preceq \pi{\tilde{\ }}} (-1)^{\ell(\pi{\tilde{\
    }}) - \ell(\tau)} \Lambda^\tau$ (where $\widebar{\tau}$ is the
complement of $\tau$ and $\pi{\tilde{\ }}$ is the conjugate of $\pi$).

It follows that, in the Lambda basis, the characteristic of the module
of parking functions is given by an alternating sum:
\begin{align*}
    \ncch(\PF_n(\seqPF)) &= \sum_{\pi \models n} \left(
                    \sum_{\tau \models \ell(\pi)} (-1)^{n - \ell(\tau)}
                        \prod_{i=1}^{\ell(\tau)}
                            \binom{\seqPFn{1 + \pi(\tau(i-1))}}{\tau_i}
                    \right) \L^\pi \,.
\end{align*}
\begin{ex}
  \begin{align*}
    \ncch(\PF_1(m^2-m+1)) &= \L^1,\qquad
    \ncch(\PF_2(m^2-m+1)) = 3\L^{11} - \L^2\\
    \ncch(\PF_3(m^2-m+1)) &= 18\L^{111} - 3\L^{1 2} - 7\L^{2 1} + \L^{3}\\
    \ncch(\PF_4(m^2-m+1)) &= 172\L^{1111} - 18\L^{112} - 36\L^{121} + 3\L^{13} -
    70\L^{211} + 7\L^{22} + 13\L^{31} - \L^{4}
  \end{align*}
\end{ex}
The coefficients once again admit a combinatorial interpretation.
\begin{prop}
  \label{prop:lambda-combi}
Let $\pi$ be a composition of $n$. The coefficient of $[\L^{\pi}]$ is the
number of non-decreasing $\seqPF$-parking functions constant on each part of
$\pi$, up to the sign.
\end{prop}
\begin{ex}
  The coefficient $[\L^{\pi}] \ncch(\PF_4(m^2-m+1))$ of the
  previous example is $7$; this is the number of non-decreasing parking functions that are
  constant on each part of the composition $22$: 
  $\boldsymbol{(} 1234 | \cdot | \cdot | \cdots \boldsymbol{)}$, 
  $\boldsymbol{(} 12 | 34 | \cdot | \cdots \boldsymbol{)}$, \ldots,
  $\boldsymbol{(} 12 | \cdot | \cdot | \cdot | \cdot | \cdot | 34 | \cdots
  \boldsymbol{)}$.
\end{ex}

%% file: incl-excl.tex
Originally we expected the formula
\begin{align*}
    \structPF(\seqPF; n) &= \sum_{\substack{\pi \models n\\\pi = \pi_1 \cdots
    \pi_k}}
    (-1)^{n-k} \binom{n}{\pi_1,\dots,\pi_k} 
    \prod_{i=1}^{k} \seqPFn{1 + \pi(i-1)}^{\pi_i}
\tag*{\cite[Theorem 4.2]{kung2003goncarov}}
\end{align*}
to be the specialization at $\L^\pi \mapsto \binom{n}{\pi_1, \cdots,
  \pi_k}$ of the non commutative characteristic of the module of
generalized parking functions.  This it turned out is not to be the case,
therefore the aim of this section is to investigate this formula, in particular
to try to find a representation theoretic interpretation of it.

First we need a few definitions;
given an ordered alphabet $A$, recall that the \emph{standardization} of a
word $w \in A^*$ is the permutation obtained by scanning iteratively $w$ from
left to right and relabeling $1, 2, \ldots, n$ the occurrences of the smallest
letters. For any non-decreasing sequence  $\seqPF$ of integers, the
$\seqPF$-\emph{standardization} of a word $w$ is the word obtained by applying
the same algorithm and then relabeling with $\seqPF(1), \seqPF(2), \cdots,
\seqPF(n)$. The word (in fact a generalized parking function) obtained is no
longer a permutation. We denote by $\stdt$ this operator.
\begin{ex}
  Let $\seqPF$ be the Catalan numbers $1, 1, 2, 5, 14, 42, \cdots$
  we have,
  \begin{align*}
   \stdt(1,4,11,1,31,1) = (1, 5, 14, 1, 42, 2).
  \end{align*}
\end{ex}
\begin{defi}
  A $\seqPF$-parking function $(Q_i)$ of size $n$ is \emph{primitive} if the
  following is verified:
  \begin{align*}
    Q_i \neq \emptyset &\iff i = \seqPF\left(1 + \sum_{j<\seqPFn{i}}
    \#Q_j\right).
  \end{align*}
  We denote by $\p_{n-1}$ the set of primitive parking
  functions of size $n$.
\end{defi}
In other words, in the sorted of $f$, all vertical paths join
$\seqPF$.

An obvious bijection between primitive parking functions and ordered
set partitions of $n$ is obtained by considering the sequence of non
empty $Q_i$ in the same order.
%
The \emph{inversion set} of a primitive parking function $f$ is :
$\inv(f) = \{(i, j) \suchthat i < j \text{ and } f(i) \geq f(j)\}$.
\begin{ex}
  Let $\seqPF$ the sequence of prime numbers, the function $(13, 2, 3, 11,
  3, 3)$ is primitive and its associated ordered set partition is $\{2\} |
  \{3,5,6\} | \{4\} | \{1\}$. The inversion set of $f$ is
  \begin{align*}
    \inv(f) &= \{(1,2), (1,3), (1,4), (1,5), (1,6), (3,5), (3,6), (4,5), (4,6),
    (5,6)\}.
  \end{align*}
\end{ex}
The collection of ordered set partitions admits a nice representation
as indexing the faces of a polytope (see~\cite{MR1311028}). In this
polytope the faces of dimension $i$ are the ordered set partitions
with $n-i$ parts. In particular the ordered partition with only one
part corresponds to the only face of dimension $n-1$.

Through the aforementioned bijection we may alternatively label the
faces of this polytope with the $\seqPF$ primitive parking functions;
the dimension of the face indexed by $f$ is then $d(f) = n - \img(f)$.

We name $\p_{n-1}$ the $n-1$ dimensional polytope of primitive parking
function of size $n$. If $\seqPF$ is strictly increasing this is the
\emph{permutohedron}.
Generalized parking functions (on $[n]$) can naturally be endowed with
the product order inherited from $\NN^n$: namely $f \leq g$ if and
only if $\forall i \leq n$, $f(i) \leq g(i)$. Seeing each face $e_f$
of $\p_{n-1}$ as the sum of parking functions lower than $f$, our main
theorem states that, by doing an inclusion-exclusion process on the
dimension of the faces, we obtain each generalized-parking function
once and only once.
\begin{thm}
  \label{thm:inc-exc}
  In the vector space $\CC\,\PF_{n}(\seqPF)$ one has
  \begin{align*}
    \sum\limits_{f \in \PF_n(\seqPF)} f
    &=
    \sum \limits_{f \in \p_{n-1}} (-1)^{d(f)} \sum \limits_{p \leq f}
    p\,.
  \end{align*}
\end{thm}
More combinatorially, one can directly count $\seqPF$-parking functions from the
previous theorem:
\begin{coro}[\expandafter{\cite[Theorem 4.2]{kung2003goncarov}}]
  \begin{align*}
    \label{formula:ky1}
    \structPF(\seqPF; n) &= \sum\limits_{\substack{Q \in \PF_n(\seqPF)\\
      Q\ \textrm{primitive}}} (-1)^{d(Q)} \prod\limits_{i=1}^{n}
      \seqPFn{i}^{\#Q_{\seqPFn{i}}},
  \end{align*}
\end{coro}
%
We prove it by using the \emph{signed involution principle}.
  The aim is, for any $\seqPF$-parking function $f$ of
size $n$, to give an involution $\invo_f$ from $\{p \in \p_{n-1}
\suchthat p \geq f\}$ into itself such that:\\[-10pt]
\begin{align*}
  \invo_f(p) = 
  \left\{
  \begin{array}{l l}
    p & \text{ if } p = \stdt(f)\\
    y\ \text{ with } d(y) = d(f) \pm 1 & \text{ otherwise.}
  \end{array}
  \right.
\end{align*}
\begin{table}
  \center \scalebox{.8}{\input{subcomplex1.tex}}
  \caption{The set $\{p \in \p_{n-1} \suchthat p \geq f\}$ for $\seqPF = \id$
and $f = (1123)$. Arrows represent a possible involution.}
\end{table}
Without loss of generality we can suppose that $f$ is non-decreasing
so that $\stdt(f) = (\seqPF(1), \cdots \seqPF(n))$. The involution is
defined implicitly from the inversion set of $p$; it is $\stdt(f)$ if
$\inv(p) = \emptyset$ and another primitive parking function of
another dimension but with same inversion set otherwise. The key
ingredient is given a primitive parking function $p$ with inversion
set $\inv(p) = I$, to understand the set
\begin{align*}
 \inv_p =
  \{q \in \p_{n-1} \suchthat \inv(q) = \inv(p) \text{ and } q \geq
p\}.
\end{align*}
For a parking function $f \in \p_{n-1}$ the dimension is an invariant of the
symmetric group action on the indices, as well as the cardinality of $\inv_p$.
We can then state the following lemma:
\begin{lemme}
  \label{lemma:distr-inv}
  let $I$ be an inversion set and $h \geq f$ a primitive $\seqPF$-parking
  function with inversion $I$ of maximum dimension. The generating series
  $G_I(f)$ of the dimensions of the faces with inversion set $I$ is
  \begin{align*}
    G_I(f) &= \sum_{\substack{\inv(p) = I\\p \geq f}} \x^{d(p)}
            = (1+\x)^{d(h)}.
  \end{align*}
\end{lemme}
%
\longversion{
\begin{proof}
  Let $h$ a primitive parking function with inversion set $I$. Thanks to the previous remark,
  we can suppose than $p$ is non-decreasing without loss of generality. Because
  $h$ is non-decreasing, the set $\inv_p$ is obtained by replacing the parts of $[n]$
  with same image by primitive parking function. Equivalently, we choose a
  subset of size $d(h)$ of indices and force $h$ to increase on those indices.
  We hence obtained $\binom{d(h)}{k}$ primitive parking functions of $\inv_h$ of
  dimension $d(h) - k$. The lemma follows if $h$ is of maximum dimension with
  inversion set $I$.
\end{proof}
}
The construction of $f_x$ is straightforward from 
Lemma~\ref{lemma:distr-inv}, which completes the proof of Theorem~\ref{thm:inc-exc}.
By adding the number of $\seqPF$-parking functions in each dimension
we get the following formula:
\begin{prop} \ \\[-20pt]
  \begin{align*}
    0 &=\sum_{k=0}^{n} (-1)^{n-k+1} \binom{n}{k} \structPF(\seqPF; k)
    \seqPF(k+1)^{n-k}.
  \end{align*}
\end{prop}
\subsection*{Frobenius characteristic investigation}
In the previous subsection the formula of Theorem 4.2
of~\cite{kung2003goncarov} is expressed combinatorially as an
alternating sum. This formula is the result of the exponential
specialization of the following non-commutative characteristic:
\begin{align*}
    \mathcal{G}(\seqPF; n) = \sum_{\pi \models n} (-1)^{n - \ell(\pi)}
    \prod_{i=1}^{\ell(\pi)} \seqPFn{1 + \pi(i-1)} \Lambda^\pi.
\end{align*}
Unfortunately this expression is not positive when expanded on the $\R$ basis.
It's therefore not the characteristic of an indecomposable $H_n(0)$-module.
Nevertheless it might still be interpretable as the characteristic of some
exact sequence of $H_n(0)$-modules.

%% file: subcomplex1.tex
\begin{tikzpicture}
  \matrix (m) [matrix of math nodes, row sep=1em, column sep=3em]{
    &&1243&&2143&&\\
    &&&1133&&&\\
    1423 & 1324 & 1234 && 2134 & 3124 & 4123\\
  };
  \path[-] (m-1-3) edge node[above] (l1143) {$1143$} (m-1-5) 
           (m-1-3) edge node[left] (l1233) {$1233$} (m-3-3)
           (m-1-5) edge node[right] (l2133) {$2133$} (m-3-5)
           (m-3-1) edge node[below] (l1323) {$1323$} (m-3-2)
           (m-3-2) edge node[below] (l1224) {$1224$} (m-3-3)
           (m-3-3) edge node[below] (l1134) {$1134$} (m-3-5)
           (m-3-5) edge node[below] (l2124) {$2124$} (m-3-6)
           (m-3-6) edge node[below] (l3123) {$3123$} (m-3-7);
  \path[<->, orange] (m-3-1) edge[bend right] (l1323)
                  (m-3-2) edge[bend right] (l1224)
                  (l1134) edge[bend right] (m-3-5)
                  (l2124) edge[bend right] (m-3-6)
                  (l3123) edge[bend right] (m-3-7)
                  (l1233) edge[bend left] (m-1-3)
                  (l1143) edge[bend left] (m-2-4)
                  (m-1-5) edge[bend left] (l2133)
                  (m-3-3) edge[loop below] (m-3-3);
\end{tikzpicture}

%% file: acknowledgement.tex
We would like to thank Vincent Pilaud and Jean Christophe Novelli for their
help and many discussions, especially about the last section. The first author
would also like to thank François Bergeron for the truly
rewarding summer internship in Montreal (supported by LIA LIRCO).